\documentclass[10pt]{amsart}
\usepackage{amscd}
\usepackage{verbatim}
\theoremstyle{plain}
\newtheorem{theorem}{Theorem}[section]
\newtheorem{lemma}[theorem]{Lemma}
\newtheorem{prop}[theorem]{Proposition}
\newtheorem{corollary}{Corollary}
\theoremstyle{definition}
\newtheorem{definition}[theorem]{Definition}
\theoremstyle{remark}
\newtheorem{remark}[theorem]{Remark}
\numberwithin{equation}{section}
\newtheorem{example}{Example}
\newcommand\ox{\mathcal O_{X}}
\newcommand\oy{\mathcal O_{Y}}

\newcommand\mo{\mathcal O}
\newcommand\lra{\longrightarrow}
\newcommand\ml{\mathcal L}
\newcommand\omtx[1][1]{\Omega^{#1}_{\tilde X}}
\newcommand\mtx{\Theta_{X}}
\newcommand\mty{\Theta_{Y}}
\newcommand\omy[1][1]{\Omega^{#1}_{Y}}
\newcommand\CC{\mathbb C}
\newcommand\PP[1][3]{\mathbb P^{#1}}
\newcommand\mn{\mathcal N}
\newcommand\lrap[1][\pi]{\xrightarrow{\rule{3mm}{0mm}#1\rule{3mm}{0mm}}}
\renewcommand\ker{\mathcal{K}\mathit{er}}
\newcommand\mt{\Theta}
\newcommand\coker{\mathcal C\mathit o\mathcal{K}\mathit{er}}
\newcommand\tyd[1][Y]{\Theta_{Y}(\log D)}
\newcommand\tyl[1][Y]{\Theta_{Y}\otimes\mathcal L^{-1}}
\DeclareMathOperator{\spec}{Spec}
\DeclareMathOperator{\codim}{codim}

\newcommand\ieq{\mathrm I_{\rm eq}}
\newcommand{\df}{H^{1}\Theta}

\begin{document}
\title[Infinitesimal deformations of double covers]
{Infinitesimal deformations of double covers of smooth algebraic varieties}
\author{S\l awomir Cynk}
\author{Duco van Straten}
\address{Instytut Matematyki\\Uniwersytetu Jagiello\'nskiego\\
  ul. Reymonta 4\\30--059 Krak\'ow\\Poland}
\curraddr{Fachbereich 6\\ Universit\"at Essen\\ D--45117 Essen\\Germany}
\email{cynk@im.uj.edu.pl}
\address{Fachbereich Mathematik 17\\ Johannes
  Gutenberg-Universit\"at\\
  Staudingerweg 9\\D--55099 Mainz\\Germany}
\email{straten@mathematik.uni-mainz.de}
\maketitle

\section{Introduction}
\label{sec:intro}
The goal of this paper is to give a method to compute
the space of infinitesimal deformations of a double cover of a smooth
algebraic variety. This research was inspired by the analysis of
Calabi--Yau manifolds that arise as smooth models of double covers of 
$\PP$ branched along singular octic surfaces (\cite{CS,Cynk}). It is
of considerable interest to determine the Hodge numbers for these manifolds,
but the methods to compute these are only available in very special cases.
For example, the results of \cite{clemens} are applicable in the case where 
the octic has only ordinary double points.
Since for a Calabi--Yau 3-fold the Hodge number $h^{1,2}$ equals the 
dimension of the space of infinitesimal deformations our approach is to study the
latter.

Let $X\lrap Y$ be a double cover of a non--singular, complete, complex 
algebraic variety $Y$
branched along a non--singular (reduced) divisor $D$. 
In the space $\df_{X}$ of all infinitesimal deformations of $X$ 
one can distinguish two subspaces
\begin{itemize}
\leftskip=8mm
\item [$T^{1}_{X\to Y}$ --]infinitesimal deformations of $X$, which are double
  covers of deformations of $Y$,
\item [$T^{1}_{X/Y}$ --]infinitesimal deformations of $X$, which are double
  covers of  $Y$.
\end{itemize}

In Proposition \ref{prop3} we give formulae for  the above two
subspaces. They have the following geometrical interpretation: 
the space $T^{1}_{X\to Y}$ is isomorphic to the space of simultaneous
deformations of $D\subset Y$, whereas $T^{1}_{X/Y}$ is isomorphic to the space of
deformations of $D$ as a subscheme of $Y$ modulo those coming from
infinitesimal automorphisms of $Y$. The space of simultaneous
deformations of $D\subset Y$ is isomorphic to the cohomology group
$H^{1}(\tyd)$, the so--called
logarithmic deformations (cf. \cite[\S~2]{kawamata})

In the space of all infinitesimal deformations of $X$ we identify
(Proposition~\ref{prop2}) a subspace 
isomorphic to $H^{1}(\mt_{Y}\otimes\ml^{-1})$, which is complementary to $T^{1}_{X\to
  Y}$, where $\ml$ is the line bundle on $Y$ defining the
double cover. We call deformations from this subspace \emph{transverse} because they
induce trivial (up to order one) deformations of the branch locus $D$.

The main result of the paper is a description of  the effect of imposing
singularities in the branch locus.
If the divisor $D$ is \emph{singular} then there exist a sequence of blow--ups
$\sigma:\tilde Y\lra Y$ and a non--singular (reduced) divisor $D^{\ast}\subset
\tilde Y$, s.t. $\tilde D\leq D^{\ast}\leq  \sigma^{\ast}D$ and $D^{\ast}$ is an even
element of the Picard group of $\tilde Y$. The double cover $\tilde X$ of  
$\tilde Y$ branched along $D^{\ast}$ is a smooth model of $X$.     

We prove that the space $T^{1}_{\tilde X\lra \tilde Y}$ can be interpreted as the space of
equisingular deformations of $D$ in deformations of $Y$,
under the additional
assumption that $Y$ is rigid it is just the space of
\emph{equisingular}  
deformations of $D$ in $Y$ (Theorem~\ref{thm1}).
Using this interpretation we give  an explicit formula for the space of
infinitesimal deformations. This formula has a particularly simple
form when the base  $Y$ is a projective space.
 In this case the equisingular
deformations can be computed from the equisingular ideal which can be 
written down in terms of the resolution of singularities
(Theorem~\ref{thm2}). The main advantage of this formula is that all computations
are carried out on $\PP[n]$ (not on the blow--up), which makes this method
very effective.

We separately study the effect of imposing singularities on the
transverse deformations. It is quite easy to compute dimension of this space,
 on the other hand their geometry may
be quite complicated. Transverse deformations of double cover induce
(second order) deformations of the branch divisor which are not
equisingular. We study some examples which exhibit the possible
phenomena. 

As a special case we study deformations of Calabi--Yau threefolds which are
non--singular models of double cover of $\PP$ branched along an octic
surface. We show that in that case the number of deformations can be computed
explicitly using computer algebra systems. This gives a method to compute the
Hodge numbers of these Calabi--Yau manifolds.
In this case dimension of the space of transverse deformations is
easily computed as the sum of genera of all curves blown--up during
the resolution. The deformations of a Calabi--Yau manifold are
unobstructed, so in that case we can
study small deformations. In this situation transverse deformations are
resolutions of deformations of double covers of $\PP$ but not double
covers of a blow--up of $\PP$ (cf. Remark~\ref{rem52}).
 
\section{Infinitesimal deformations}
\label{sec:id}
An \emph{infinitesimal deformation} of $X$ is any scheme $X'$ flat over the
ring of dual numbers $\mathbb D=\CC[t]/[t^{2}]$ such that
$X'\otimes_{\mathbb D}\CC\cong X$.
If the variety $X$ is smooth then the space of infinitesimal deformations 
is isomorphic to the cohomology group 
$H^{1}\mtx$ of the tangent bundle $\mtx$.

Let $\pi:X\lra Y$ be a double cover of a \emph{smooth} algebraic variety branched
along a smooth divisor $D$. 
The cover $\pi$ is not determined by $D$ itself, we have also to fix a line
bundle $\ml$ on $Y$ s.t. $\pi_{\ast}\ox\cong\oy\oplus\ml^{-1}$. This $\ml$
satisfies $\ml^{\otimes2}\cong\oy(D)$.
Since the map $\pi$ is finite we have
\(H^{i}(\mtx)\cong H^{i}(\pi_{\ast}\mtx)\). From \cite[Lem.~3.16]{EV} we get
$\pi_{\ast}\mtx\cong\mty\otimes\ml^{-1}\oplus\tyd$, where $\tyd$ is the
 sheaf of logarithmic vector fields which is 
defined by the following exact sequence
The sheaf $\tyd$ is the kernel of the natural  restriction map
\(\mt_{Y}
\lra\mn_{D|Y}\)
and so it is the subsheaf of the tangent bundle $\mty$ consisting of those
vector fields which carry the ideal sheaf of $D$ into itself.
\begin{equation}
  \label{cd1}
0\lra
\tyd\lra\mt_{Y}
\lra\mn_{D|Y}\lra 0.  
\end{equation}

This gives immediately
\begin{prop}
\label{prop2}
  \[H^{1}\mtx\cong H^{1}(\tyd)\oplus H^{1}(\mty\otimes\ml^{-1}). \]
\end{prop}

\begin{prop}\rule{0cm}{0cm}\samepage  
  \label{prop3}
  \begin{itemize}
\item [(a)] \(H^{1}(\tyd)\cong 
\coker(H^{0}\mty\lra H^{0}\mn_{D|Y})
  \oplus\ker(H^{1}\mty\lra H^{1}\mn_{D|Y}),\)\medskip
  \item [(b)]   $H^{1}(\tyd)$ is isomorphic to the space $T^{1}_{X\to Y}$ of
  infinitesimal deformations of $X$ which are double
  covers of deformations of~$Y$,\medskip
\item [(c)] $\coker(H^{0}\mty\lra H^{0}\mn_{D|Y})$ is isomorphic to the space
  $T^{1}_{X/Y}$ of infinitesimal deformations of $X$ which are double
  covers of ~$Y$.
  \end{itemize}

\end{prop}

\begin{proof}
The cohomology exact sequence derived from \eqref{cd1} yields
\[H^{0}\mt_{Y}\lra H^{0}\mn_{D|Y}\lra\df_{Y}(\log D)\lra
H^{1}\mt_{Y}\lra  H^{1}\mn_{D|Y}\]
which proves (a).

The maps
\(H^{0}\mt_{Y}\lra H^{0}\mn_{D|Y}\) and
\(H^{1}\mt_{Y}\lra  H^{1}\mn_{D|Y}\) have quite obvious interpretations. 
The first one associate to an infinitesimal automorphism of $Y$ an
infinitesimal deformation of $D$ in $Y$, consequently
\(\coker (H^{0}\mt_{Y}\lra H^{0}\mn_{D|Y})\)
is the space of deformations of $D$ as a subscheme of  $Y$ modulo
automorphisms of $Y$.
The second map 
\(H^{1}\mt_{Y}\lra  H^{1}\mn_{D|Y}\)
gives for a deformation of $Y$ the obstruction to lift it to a deformation of
$D$. Indeed, from the diagram
\[
\begin{CD}
@. {H^{1}\Theta_{Y}} \\
@. @VVV\\
H^{1}\Theta_{D}@>>> H^{1}(\Theta_{Y}\otimes\mo_{D})@>>> H^{1}\mn_{D|Y}
\end{CD}
\]
we see that if an element of $H^{1}\mty$ belongs to 
$\ker(H^{1}\mty\lra H^{1}\mn_{D|Y})$ then its image in $H^{1}(\Theta_{Y}\otimes\mo_{D})$
lies in the image of $H^{1}\Theta_{D}$.
Consequently $H^{1}(\tyd)$ is isomorphic to the space of simultaneous
deformations of $D\subset Y$, i.e.\ pairs $D'\subset Y'$ such that $D'$ is an
infinitesimal deformation of $D$ and  $Y'$ is an
infinitesimal deformation of $Y$.
Let $X'$ be an infinitesimal deformation of $X$ which is  a double
cover of a deformation $Y'$ of $Y$. Denote by $D'\subset Y'$ the branch
locus. Restricting to the central fiber we find that
$D'\otimes_{\mathbb D}\CC\cong D$ and so $D'$ is an infinitesimal
deformation of $D$. 

Conversely if $D'\subset Y'$ are deformations of $D\subset Y$ the $D'$ is even
and there exists a unique line bundle $\ml'$ on $Y'$ such tht $\ml'|Y\cong\ml$
and $(\ml')^{\otimes2}\cong\mo_{Y'}(D')$, The line bundle $\ml'$ is defined by the square root
of the transition functions of an extension to $X'$ of
$\ml^{\otimes2}$. The transition functions of $\ml'$ are of the form
$f^{2}+\epsilon g$, and the square root  equals $f+\frac12\epsilon g$.
The line bundle $\ml'$ defines a
double cover $X'\lra Y'$ branched along $D'$, restricting to the
central fiber we find that $X'\otimes_{\mathbb D}\CC$ is a double
cover of $Y$ branched along $D$ defined by the line bundle
$\ml$. This means that 
  $X'\otimes_{\mathbb D}\CC$ is isomorphic to $X$ and so $X'$ is a
  deformation of $X$. This proves (b), and also (c) easily follows.
\end{proof}
\begin{corollary}\rule{0cm}{0cm}
  \begin{itemize}
  \item [(a)] Every deformation of $X$ is a double cover of a deformation of $Y$ iff
    $H^{1}(\mty\otimes\ml^{-1})=0$. 
  \item [(b)] Every deformation of $X$ is a double cover of $Y$ iff
    $H^{1}(\mty\otimes\ml^{-1})=0$ and the map $H^{1}\mty\lra H^{1}\mn_{D|Y}$
    is injective (e.g. $Y$ is rigid).
  \end{itemize}
\end{corollary}
\begin{remark}
  $H^{1}(\mty\otimes\ml^{-1})$ is isomorphic to the space of infinitesimal
  extensions of $Y$ by $\ml^{-1}$. If $(Y',\mathcal F)$ is any such extension, then 
$\spec(\mo_{Y'}\oplus\mathcal F)$ is an infinitesimal deformation of $X\cong
  \spec(\mo_{Y}\oplus\ml^{-1})$. At the beginning of the section
  \ref{sec:sym} we give a more geometric interpretation.
 
\end{remark}
\begin{example}
  Let $Y=\PP[n]$ ($n\geq2$) and let $D$ be a smooth hypersurface of degree
  $2d$. Then $H^{1}\mty=0$ for any $n,d $ and $H^{1}(\mty\otimes\ml^{-1})=0$ with
  the only exception of $d=6,\;n=2$. So for $(d,n)\not=(6,2)$ every
  infinitesimal deformation of a double cover of $\PP[n]$ branched along a
  degree $d$ \emph{smooth} hypersurface is again a double cover of $\PP[n]$ branched along a
  smooth hypersurface of the same degree.

In the case $d=6,\;n=2$ the dimension of the space of infinitesimal
deformations of a K3 surface is 20, whereas the dimension of the family of double sextic K3
surfaces equals $\dim T^{1}_{X\to Y}=\dim T^{1}_{X/Y}=19$. 
\end{example}

\begin{example}
\label{ex:twosurf}
  Let $D_{1}$ and $D_{2}$ be two surfaces if $\PP[3]$ of degree $d_{1}$ and
  $d_{2}$ intersecting transversely along a smooth curve $C$. Let
  $Y=\text{Bl}_{C}\PP$ be the blow-up of $\PP$ along $C$,
  $D^{\ast}=\tilde{D_{1}}+\tilde{D_{2}}$, where $\tilde{D_{i}}$ is the strict
  transform of $D_{i}$. Consider the double cover $\pi:X\lra Y$ of $Y$
  branched along $D^{\ast}$. The exceptional divisor of the blow--up $E$ and its
  pullback to the double cover $E_{1}$ are ruled surfaces over $C$. Simple
  computations yields
\(h^{1}(\Theta_{Y}(\log D^{\ast}))=\binom{d_{1}+3}3+\binom{d_{2}+3}3-17\)
and \(h^{1}(\mty\otimes\ml^{-1})=h^{0}\mo_{C}\left(\tfrac12\left({d_{1}
      +d_{2}}\right)\right)\).

We can give an explicit description of deformations of the double cover which
are not double cover of a deformation of $Y$. Namely
$H^{0}\mo_{C}\left(\tfrac12\left({d_{1}+d_{2}}\right)\right)$ is the space of restrictions to
$C$ of degree $\frac12(d_{1}+d_{2})$ surfaces in $\PP$. Generic such a surface
gives a deformation of the surface $D_{1}+D_{2}$ which replace the double
curve by $\frac12d_{1}d_{2}(d_{1}+d_{2})$ nodes. The double cover of $\PP$
has a double curve along  which it is locally isomorphic to a product of a
node and a line ($cA_{1}$ singularity) which deforms to a set of nodes. This family
admits a simultaneous resolution which replace the ruled surface $E_{1}$ by 
 $\frac12d_{1}d_{2}(d_{1}+d_{2})$ lines (for the special element of the family
 the resolution is a blow--up of a double curve whereas for the general one it
 is a small resolution of nodes).

\end{example}

\section{Resolution of singularities of a double cover}
\label{sec:sing}

Let $Y$ be a non--singular complex algebraic variety and let
$D$ be a divisor on $Y$ 
which is even as an element of the Picard
group. 
Let $\ml$ be a line bundle on $Y$ s.t. $\ml^{\otimes2}\cong\mo_{Y}(D)$.
Consider $X\lrap Y$ the
double cover of $Y$ branched along  $D$ defined by $\ml$.
$X$ is non--singular iff $D$ is non--singular, if $D$ is
singular then the singularities of $X$  are in one--to--one
correspondence with singularities of $D$.

We can resolve singularities of $X$ by a special resolution of $D$. 
For any birational morphism $\sigma:\tilde Y\lra Y$ we have
\((\sigma^{\ast}D)=\tilde D+\sum_{j}m_{j}E_{j}\)
(where $\tilde D$ is the strict transform of $D$, 
$E_{j}$ are the $\sigma$--exceptional divisors and
$m_{j}\geq0$). Therefore the divisor
\[D^{\ast}=\tilde D+\sum_{2\not\;\vert m_{j}}E_{j}=
\sigma^{\ast}D-2\sum_{j}\left\lfloor\frac {m_{j}}2\right\rfloor E_{j}\]
is reduced and even. In fact it is the only reduced and even   divisor
satisfying
\[\tilde D\leq D^{\ast}\leq \sigma^{\ast}D.\]

Let $\tilde X\lrap[\tilde \pi]\tilde Y$ be the double cover branched along
$D^{\ast}$ defined by
$\sigma^{\ast}\ml\otimes\mo_{Y}(-\sum_{j}\lfloor\frac {m_{j}}2\rfloor E_{j})$, 
we can find a birational morphism
$\tilde X\lrap[\rho]X$ which fits into the following commutative diagram
\[\begin{CD}
  \tilde X@>\rho>>X\\
@V\tilde\pi VV@VV\pi V\\
\tilde Y@>\sigma>>Y
\end{CD}\]

From the Hironaka desingularization theorem it follows that we can find a sequence
of blow--ups with smooth centers $\sigma:\tilde Y\lra Y$ such that $D^\ast$ is a
smooth divisor, and consequently which gives a resolution of singularities of
the double cover.

\section{Equisingular deformations}
\label{sec:ed}

Let $\sigma:\tilde Y\lra Y$ be a resolution of singularities of $D$ as
explained in the previous section. 
In this section we shall study the infinitesimal deformations from
$H^{1}(\mt_{\tilde Y}(\log D^{\ast}))$. 

Before going to the general case consider first a single blow--up
\(\sigma:\tilde Y\lra Y\)
of a smooth subvariety $C\subset Y$, denote by $E$ the exceptional locus of
$\sigma$.

Using the Leray spectral sequence we compute that
\begin{eqnarray*}
  H^{0}\Theta_{\tilde Y}&\cong&\ker(H^{0}\Theta _{Y}\lra H^{0}\mn_{C|Y})\\
  H^{1}\Theta_{\tilde Y}&\cong&\coker(H^{0}\Theta _{Y}\lra H^{0}\mn_{C|Y})\oplus
  \ker(H^{1}\Theta _{Y}\lra H^{1}\mn_{C|Y})
\end{eqnarray*}
The above formulas have  nice geometric interpretations.
The space of infinitesimal automorphisms
of $\tilde Y$ is isomorphic to the space of 
infinitesimal automorphisms of $Y$ which fix the subvariety $C$.
The space of infinitesimal deformations of $\tilde Y$
 is isomorphic to the direct sum of the 
space infinitesimal deformations of $C$ as a subscheme of $Y$ modulo
those coming from infinitesimal automorphisms of $Y$ and the space of
infinitesimal deformations of $Y$ which can be lifted to a deformation
of $C\subset Y$. 
This vector space is isomorphic to  the space of simultaneous deformations of $C\subset Y$
modulo those coming from infinitesimal automorphisms of $Y$ (the first
summand controls the deformations of $C$, while the second one -- the
deformations of $Y$).

Recall that  $D^{\ast}=\sigma^\ast D-mE$ and so  $H^{0}\mn_{D^{\ast}|\tilde Y}$ is isomorphic to the
subspace of $H^{0}\mn_{D|Y}$ corresponding to those infinitesimal deformations of
$D$ in $Y$ that have multiplicity at least $m$ along $C$. Consequently the cokernel 
$\coker(H^{0}\mt_{\tilde Y}\lra H^{0}\mn_{D^{\ast}|\tilde Y})$ is the space of these infinitesimal
deformations modulo infinitesimal automorphisms of $Y$ which fix $C$.
In a similar manner the kernel 
$\ker(H^{1}\mt_{\tilde Y}\lra H^{1}\mn_{D^{\ast}|\tilde Y})$ is the space of simultaneous
deformations of $C\subset Y$ which can be extended to a simultaneous deformation
$C'\subset D'\subset Y'$ of $C\subset D\subset Y$ such that $D'$ has
multiplicity at least $m$ along $C'$.
As in the formula for $H^{1}(\Theta_{\tilde Y})$ the above two 
subspaces gives the space of all simultaneous deformations
$D'\subset Y'$ of $D\subset Y$, which can be extended to a deformation
$C'\subset D'\subset Y'$ of $C\subset D\subset Y$ such that $D'$ has 
multiplicity at least $m$ along $C'$.

Going back to the general case, let $\sigma:\tilde Y\lra Y$ be a sequence
$\sigma=\sigma_{n-1}\circ\dots\circ\sigma_{0}$, where $\sigma_{i}:Y_{i+1}\lra
Y_{i}$ is a blow--up of a smooth subvariety $C_{i}\subset Y_{i}$ such that 
$D^{\ast}$ is smooth, $Y_{0}=Y$, $Y_{n}=\tilde Y$. Let $m_{i}$ be an integer such that
$D_{i+1}^{\ast}=\sigma_{i}^{\ast}D_{i}^{\ast}-m_{i}E_{i}$, where $E_{i}\subset
Y_{i+1}$ is the exceptional divisor of $\sigma_{i}$.
Applying the above description to every $\sigma_{i}$ separately we see
that any deformation from $H^{1}(\mt_{\tilde Y}(\log D^{\ast}))$ gives
(inductively) a deformation of $D^{\ast}_{i}\subset Y_{i}$. From the
above description we conclude the following 
\begin{theorem}
  \label{thm1}
$H^{1}(\mt_{\tilde Y}(\log D^{\ast}))$ is isomorphic to the space of
simultaneous deformations of $D\subset Y$ which have simultaneous resolution
i.e.\ which can be lifted to deformations of $C_{i}\subset
D^{\ast}_{i}\subset Y_{i}$ in such a way that
the multiplicity of the deformation of $D_{i}^{\ast}$ along
deformation of $C_{i}$ is at least $m_{i}$.
\end{theorem}
\begin{definition}
  We call an infinitesimal deformation of $D$ in $Y$  \emph{equisingular} if
  it satisfies the assertion of the above theorem.
\end{definition}
The above theorem is particularly useful  when we have an explicit description of
infinitesimal deformations of $Y$, for instance  when $Y$ is rigid.
\begin{corollary}
  If the variety $Y$ is rigid then the space of 
  equisingular deformations of $D$ in $Y$  is isomorphic to 
  $H^{1}(\mt_{\tilde Y}(\log D^{\ast}))$.
\end{corollary}
\begin{remark}
  The notion of an equisingular deformation is relative to a fixed resolution
  of singularities. 

The notion of equisingularity is relative to a fixed embedded resolution
of singularities and is equivalent to the existence of simultaneous
resolution,  which is the definition formulated by Wahl in \cite{Wahl} and 
Kawamata in \cite{kawamata}.

\end{remark}
\begin{example}
  Let $X\subset \PP[N]$ be a hypersurface with a cusp ($A_{2}$
  singularity). Then $X$ has two natural resolutions: the minimal (one blow--up
  of the double point) and the log--resolution (two blow--ups: first the double
  point and then the intersection of the exceptional divisor with the strict
  transform). 
These two resolutions lead to different spaces of equisingular
deformations. For the first one equisingular are deformations with a double
point whereas for the second one -- with a cusp. The explanation is that 
the only information that we can get from first resolution (minimal)  is that
we have a hypersurface with a double point. From the second we know that the
strict transform is tangent to the exceptional locus which means that the
second derivative of the equation vanish along a line.
\end{example}

If $Y$ is  rigid we can use Theorem~\ref{thm1} to 
give a more direct description of the
space of equisingular deformations of $D$ in $Y$.
Consider an embeded resolution $\sigma:\tilde Y\lra Y $ 
of $D$ in $Y$ by a sequence of
blow--ups with smooth centers. More precisely assume that
$\sigma=\sigma_{n-1}\circ\dots\circ\sigma_{0}$, where $\sigma_{i}:Y_{i+1}\lra
Y_{i}$ is a blow--up of a smooth subvariety $C_{i}\subset Y_{i}$,
$Y_{0}=Y$, $Y_{n}=\tilde Y$. Denote by $E_{i}\subset Y_{i+1}$ the
exceptional divisor of $\sigma_{i}$, let $m_{i}$ be a nonnegative
integer such that
$D_{i+1}^{\ast}=\sigma_{i}^{\ast}D_{i}^{\ast}-m_{i}E_{i}$, where
$D_{i}^{\ast}$ is an effective divisor in $Y_{i}$ and
$D_{0}^{\ast}=D$. Assume that the divisor $D^{\ast}:=D^{\ast}_{n}$
is non--singular. 

Let $\mathcal I(C_{i})$ be the ideal sheaf of of $C_{i}$ in $Y_{i}$
and  let $\tilde
{\mathcal I}_{i}^{m_{i}}$ denotes (for nonnegative integer $m_{i}$)
the push--forward
$(\sigma_{i-1}\circ\dots\circ\sigma_{0})_{\ast}(\mathcal
I(C_{i})^{m_{i}})$ to $Y$ of the $m_{i}$--th power of $\mathcal
I(C_{i})$.
Denote by $\mathcal J_{i}$ the image of the homomorphism 
$\Theta_{Y_{i}}\otimes\mo_{D^{\ast}_{i}}\lra\mn_{D_{i}^{\ast}|Y_{i}}$
and by $\tilde {\mathcal J}_{i}$ its pushforward
$(\sigma_{i-1}\circ\dots\circ\sigma_{0})_{\ast}(\mathcal{J}_{i})$ to $Y$.
Let $\mathbf{J}$ denote the image of the map $H^{0}(\Theta _{Y})\lra
H^{0}\mn_{D|Y}$ induced by the exact sequence~\eqref{cd1}.

\begin{theorem}
  \[H^{1}(\mt_{\tilde Y}(\log D^{\ast}))\cong 
  \bigcap_{i=0}^{n-1}\Biggl( H^{0}
    \biggl(\left(\tilde{\mathcal I}^{m_{i}}_{i}\otimes\mn_{D|Y}\right)+
    \tilde{\mathcal J}_{i}\biggr)\Biggr)\bigg/\mathbf{J}.\]
\end{theorem}
\begin{proof}
By theorem~\ref{thm1} $H^{1}(\mt_{\tilde Y}(\log D^{\ast}))$ 
is the space of equisingular deformations of $D$ in $Y$.
As the equisingularity condition can be verified for each center of
blow--up separately we shall study  a single blow--up 
$\sigma_{i}:Y_{i+1}\lra Y_{i}$.
We have to find deformations of $D^{\ast}_{i}$ which vanish to order
$m_{i}$ along a deformation of $C_{i}$. Since every deformation of
$C_{i}$ is locally given by a vector field $v$ we can take a
deformation of $D^{\ast}_{i}$, transform it back by $-v$ and verify if
the result vanishes along $C_{i}$ to order $m_{i}$.

Equivalently  we can start with an
  infinitesimal deformation of $D^{*}_{i}$ vanishing to order 
  $m_{i}$ along $C_{i}$ and transform it by
  $v$. Locally this deformation is given by $f+\epsilon g$, where $f$ is a
  local equation of $D^{*}_{i}$ and  $g\in \mathcal I(C_{i})^{m_{i}}$. So we have to
  substitute $x+\epsilon\cdot v(x)$ in $f+\epsilon  g$. Taking into account
  $\epsilon^{2}=0$ we get \[f(x)+\epsilon\left(f'(x)\cdot v(x)+g(x)\right)\]
  so the deformation is given by the element $f'(x)\cdot v(x)+g(x)$ 
  of $\mathcal I(C_{i})^{m_{i}}+\mathcal J_{i}$.
  Pushing--forward the above formula to $Y$ proves that the space of
  equisingular deformations of $D$ in $Y$ is isomorphic to 
  $\bigcap\limits_{i=0}^{n-1}\biggl( H^{0}
    \biggl(\left(\tilde{\mathcal I}^{m_{i}}_{i}\otimes\mn_{D|Y}\right)+
    \tilde{\mathcal J}_{i}\biggr)\biggr)$. To get 
    $H^{1}(\mt_{\tilde Y}(\log D^{\ast}))$ we have to mod out by the
    space of deformations of $D$ induced by infinitesimal
    automorphisms of $Y$ i.e $\mathbf{J}$. 
\end{proof}

We shall study in more detail the 
case when $Y=\PP[N]$, in this situation every coherent sheaf on $Y$ is
given by a graded module over $\CC[X_{0},\dots,X_{N}]$ which makes
computations much simpler. Observe first that  $\mathbf{J}$ 
equals $\left(\mathrm J_{F}\right)_{d}$ the space of degree
$d$ forms in the Jacobian ideal $\mathrm J_{F}:=(\frac{\partial
    F}{\partial Z_{0}},\dots,\frac{\partial
    F}{\partial Z_{N}})$ of $F$, where $d$ is the degree of $D$ and
  $F$ is its homogeneous equation.
If  $C_{i}$ is not contained in the exceptional locus of
$\sigma_{i-1}\circ\dots\circ\sigma_{0}$ then $\tilde{\mathcal I}_{i}^{m_{i}}$
  equals the symbolic power ${\mathcal I}(\tilde
  C_{i})^{(m_{i})}$ of the ideal sheaf ${\mathcal I}(\tilde
  C_{i})$ of $\tilde C_{i}$, where $\tilde C_{i}$ is the image of $C_{i}$ in
  $  \PP[N]$. The ideal sheaf ${\mathcal I}(\tilde
  C_{i})$ is the sheaf associated to the homogeneous ideal
  $\mathrm{I}(\tilde C_{i})$ of $\tilde C_{i}$. Let $\rm J^{i}_{F}$ be
  the $\CC\left[X_{0},\dots,X_{N}\right]$--module associated to the
  sheaf $\tilde{\mathcal J}_{i}$ (in fact $\rm J^{i}_{F}$ is an ideal
  in $\CC[X_{0},\dots,X_{N}]$ containing $\mathrm J_{F}$).
  Define the \emph{equisingular ideal} of $D$ in $\PP[N]$ (w.r.t. $\sigma$) by
  \[\ieq(D)=\bigcap_{i=0}^{n-1}
  \left(\rm I(\tilde C_{i})^{(m_{i})}+J_{F}^{i}\right).\]

\begin{theorem}
\label{thm2}
The space of equisingular deformations of $D$ is isomorphic to the space of
degree $d$ forms in the quotient of the equisingular ideal modulo 
the Jacobian ideal
\[H^{1}(\mt_{\tilde Y}(\log D^{\ast}))\cong \left(\text {\rm
    I}_{eq}(D)/\mathrm J_{F}\right)_{d}.\] 
\end{theorem}

  If $C_{i}$ is contained in the exceptional locus of
  $\sigma_{i-1}\circ\dots\circ\sigma_{0}$  then
  the points of $C_{i}$ do not correspond to ordinary points of $Y$ but to
  infinitely near points. Vanishing at an infinitely near point gives
  on $\PP[n]$ some tangency condition, which has to be computed in
  local coordinates (cf.\ example~\ref{ex14}). 

The above theorem represents the space of equisingular 
deformations as a quotient of
two subspaces of the space of degree $d$ homogeneous forms, and so it gives a
very effective tool for explicit computations. It is particularly suitable for
computations with computer algebra systems.
The main difficulty is to find the ideal $\rm J^{i}_{F}$, we have
to describe the vector fields on $Y_{i}$, unfortunately they
push-down to rational vector fields on $\PP[N]$. If we are able to find those
rational vector fields on $\PP[N]$ which lift to regular vector fields on
$Y_{i}$ then we can compute the ideal $\mathrm J_{F}^{i}$ in
$\CC[X_{0},\dots,X_{N}]$ which contains regular functions generated by
results of applying those vector fields to the equation of $D$.  
If we have an explicit descriptions of the resolution $\sigma$ we can
use local coordinates to compute $\rm J^{i}_{F}$.
Consider the map $\sigma^{(i)}:Y_{i}\lra Y$. The regular vector fields
are transformed to $Y$ by applying differential of $\sigma^{(i)}$. The
results on $Y$ are locally given by rows of the jacobian matrix of
$\sigma^{(i)}$ pushed to $Y$. The same can be done by lifting the
equation of $D$ to $Y^{i}$, taking partial derivatives and pushing to $Y$.

In many cases the rational vector fields do not appear, which makes the
computation of the equisingular ideal much simpler. 
Let $D=\bigcup\limits_{i}D_{i}\subset \PP$ be an
arrangement as defined in \cite {CS}, i.e.\ a sum of smooth components
$D_{i}$ such that the components $D_{i}$ and $D_{j}$ (for $i\not=j$)
intersect transversally along a smooth curve $C_{ij}$, the curves $C_{ij}$
and $C_{kl}$ are either equal, or disjoint, or they intersect transversely
(locally $D$ looks like a plane arrangement). Let $\sigma$ be a
natural resolution of $D$.
Let $\sigma:\tilde Y\lra Y$ be the following resolution of $D$, 
first we blow--up the $p$--fold points that do not lie on a 
$p-1$--fold curve, then the multiple curves. Denote by
$C_{i}\;(i=0,\dots,n-1)$ the multiple points and curves of the
arrangement and by $m_{i}$ the corresponding multiplicities.

\begin{lemma}
   \[\ieq(D)=\bigcap_{i=0}^{n-1}
  \left(\rm I(C_{i})^{m_{i}}+J_{F}\right).\]
\end{lemma}
\begin{proof}
  Equisingular deformations are given by 
  arrangements of the same combinatorial type. They are given by
  deformations of the components of
  $D$ which preserves the multiplicities at points and curves. Clearly
  the centers of successive blow--ups can be interpreted as subsets of
  $\PP$, and the deformations of the centers can be obtained as the
  deformations in $\PP$. The same can be easily computed in local coordinates.
\end{proof}

\section{Transverse deformations}
\label{sec:sym}

Proposition~\ref{prop2} gives a decomposition of the group of deformations of a
double cover into two subgroups. In previous section we studied the first
summand (containing those deformations which are double covers), 
now we shall concentrate on the second one isomorphic to $H^{1}(\tyl)$.
We shall call deformations from this subspace \emph{transverse}.
A double cover of $Y$ branched along a divisor $D$ can be given as a
hypersurface $t^{2}=s$ in the total space of the line bundle
$\mo_{Y}(D)$, where $s$ is a section defining $D$. Transverse
deformations of a double cover corresponds to the deformations of
the type $t^{2}+2\epsilon tf=s$. Those deformations do not give
non--trivial first order deformations of  $D$, as we can write the
deformation locally as $(t+\epsilon f)^{2}=s+(\epsilon f)^{2}$, which is zero
because $\epsilon ^{2}=0$. On the other hand, one can use this to formally
\emph{represent} such transverse deformations as particular
second--order deformations of $D$, which is very useful in practice.
This also explain the name transverse.

Since we have the following  isomorphisms
\[H^{1}(\tyl)\cong H^{1}(\omy[n-1]\otimes K_{Y}^{\lor}\otimes\ml^{-1})\cong 
\left(H^{n-1}(\omy\otimes K_{Y}\otimes\ml)\right)^{\lor}\]
in many cases (for instance\ when $D$ is a smooth divisor in a weighted projective
space) it is easy to compute to dimension of this vector space.
We shall study the effect on $h^{1}(\tyl)$ of introducing
singularities in the branch locus of $D$, so consider a  blow--up
\(\sigma:\tilde Y\lra Y\)
of a smooth subvariety $C\subset Y$, denote by $E$ the exceptional locus of
$\sigma$, and let $m$ be such that
$D^{\ast}=\sigma^{\ast}D-mE$. Since $m=2\left\lfloor\frac{\text{\rm
      mult}_{D|C}}2\right\rfloor$, it is an even number and define
$\tilde\ml:=\sigma^{\ast}\ml\otimes\mo_{\tilde Y}(-\frac m2E)$. $\tilde\ml$ is the
line bundle on $\tilde Y$ defining the double cover, so our goal is to compare 
$h^{1}(\tyl)$ with $h^{1}(\mt_{\tilde Y}\otimes \tilde\ml^{-1})$

We have the following exact sequence
\begin{eqnarray*}
  0\lra\sigma^{\ast}(\Omega^{1}_{Y}\otimes\ml\otimes
K_{Y})\otimes \mo_{\tilde Y}(kE)\lra \Omega^{1}_{\tilde Y}\otimes\tilde\ml\otimes
K_{\tilde Y}\lra\\
\lra\Omega^{1}_{E/C}\otimes\mo_{E}(-k)\otimes \sigma^{\ast}(\ml\otimes
K_{Y})\lra0, 
\end{eqnarray*}
where $k=$codim$_{Y}C-\frac m2-1$.
Now, we can use the Leray spectral sequence to compute the required
cohomologies, the resulting formulas depends on the actual value of $k$.

The most complicated is the  case when \underline{$k<0$}. Although in this case 
\[R^{i}\sigma_{\ast}(\Omega^{1}_{\tilde Y}\otimes\tilde\ml\otimes
K_{\tilde Y})=0,\qquad\text{for }\quad i>0\]
but on the other hand the direct image 
\[\sigma_{\ast}(\Omega^{1}_{\tilde Y}\otimes\tilde\ml\otimes K_{\tilde Y})\]
is not locally free, so we cannot say too much in that case..

The easiest is the case when \underline{$k>0$}. Since $\codim_{Y}C>k+1$ simple
computations show that
\begin{eqnarray*}
  &&\sigma_{\ast}(\mo_{\tilde Y}(kE))=\oy,R^{i}\sigma_{\ast}(\mo_{\tilde Y}(kE))=0
  \quad\text{for }i\geq1\\
  &&R^{i}\sigma_{\ast}(\Omega^{1}_{E/C}(-k))=0,\quad\text{for }i\geq0.
\end{eqnarray*}
Using the projection formula we get 
\begin{eqnarray*}
  &&\sigma_{\ast}(\Omega^{1}_{\tilde Y}\otimes\tilde\ml\otimes K_{\tilde Y})\cong
  \omy\otimes\ml\otimes K_{Y},\\
  &&R^{i}\sigma_{\ast}(\Omega^{1}_{\tilde Y}\otimes\tilde\ml\otimes K_{\tilde Y})=0\quad\text{for
  }i\geq1
\end{eqnarray*}
and by the Leray spectral sequence 
\[H^{1}(\mt_{\tilde Y}\otimes\tilde \ml^{-1})\cong H^{1}(\tyl).\]
\medskip
 
The most interesting is the case when \underline{$k=0$} (crepant resolution). Since 
\begin{eqnarray*}
  &&\sigma_{\ast}(\mo_{\tilde Y})=\oy,
  R^{i}\sigma_{\ast}(\mo_{\tilde Y})=0\quad\text{for}\quad i\geq1,\\
  &&\sigma_{\ast}(\Omega^{1}_{E/C})=0,
  R^{1}\sigma_{\ast}(\Omega^{1}_{E/C})\cong \mo_{C},
  R^{i}\sigma_{\ast}(\Omega^{1}_{E/C})=0\quad\text{for}\quad i\geq2
\end{eqnarray*}
applying the projection formula and using the above exact sequence we get
\begin{eqnarray*}
  &&\sigma_{\ast}(\Omega^{1}_{\tilde Y}\otimes\tilde\ml\otimes K_{\tilde Y})\cong
  \omy\otimes\ml\otimes K_{Y},\quad
  R^{1}\sigma_{\ast}(\Omega^{1}_{\tilde Y}\otimes\tilde\ml\otimes
  K_{\tilde Y})\cong\mo_{C}\otimes\ml\otimes K_{Y},\\
  &&R^{i}\sigma_{\ast}(\Omega^{1}_{\tilde Y}\otimes\tilde\ml\otimes K_{\tilde Y})=0\quad\text{for
  }i\geq2.
\end{eqnarray*}
Now, by the Leray spectral sequence
\[H^{n-1}(\Omega^{1}_{\tilde Y}\otimes\tilde\ml\otimes K_{\tilde Y})\cong
H^{n-1}(\omy\otimes\ml\otimes K_{Y})\oplus H^{n-2} (\mo_{C}\otimes\ml\otimes
K_{Y})\]
and by Serre duality
\begin{eqnarray*}
&&H^{1}(\mt_{\tilde Y}\otimes\tilde \ml^{-1})\cong H^{1}(\tyl)
\qquad\text{if}\quad \codim_{Y}C<n-2,\\
&&H^{1}(\mt_{\tilde Y}\otimes\tilde \ml^{-1})\cong H^{1}(\tyl)\oplus
H^{0}(\det\mn_{C}\otimes\ml^{-1}) \\
&&\rule{6.7cm}{0cm}\text{if}\quad \codim_{Y}C=n-2.
\end{eqnarray*}
As a special case we get the following Proposition (notations are as before
Theorem~\ref{thm1})
\begin{prop}
  If $K_{Y}=\ml^{-1}$ and $\sigma:\tilde Y\lra Y$ is a sequence of blow--ups satisfying 
$\frac 12D^{\ast}+K_{\tilde Y}=\sigma^{\ast}(\frac 12D+K_{Y})$
then
\[h^{1}(\mt_{\tilde Y}\otimes \tilde\ml^{-1})=h^{1}(\tyl)+\sum_{\codim C_{i}=2}h^{0}(K_{C_{i}}).\]
\end{prop}
Observe that in the latter case $m=2$, which means that we are considering a
blow--up of a subvariety of codimension 2 such that the multiplicity of the
divisor along it is 2 or 3.  We shall give a geometric description of
transverse deformations in that case.

Assume first that the multiplicity of $D$ along $C$ is 2. 
If $D=D_{1}+D_{2}$ is a sum of two smooth divisors
intersecting transversely along $C$ we have
$H^{0}(\det\mn_{C}\otimes\ml^{-1})\cong H^{0}\ml$. Let $f_{i}\in
H^{0}\mo_{Y}(D_{i})$  be a section defining $D_{i}$.
For any section $f\in H^{0}\ml$ we consider the divisor $D_{\epsilon}=\{f_{1}\cdot
f_{2}+(\epsilon f)^{2}=0\}$. The family of double covers of $Y$ branched
along  $D_{\epsilon}$ have
simultaneous resolution of singularities. If the section $f$ does not vanish
along $C$ then the divisor $D_{\epsilon}$ does not contain $C$. If the zero set
of $f$ intersects $C$ transversely, then the singular locus of $D_{\epsilon}$
has codimension 2 in $D_{\epsilon}$, $D$ has ``compound nodes'' at
$\{f_{1}=f_{2}=f=0\}$. 
Moreover $D_{\epsilon}$ is irreducible and admits
a small resolution. 

If  a component of the  double points locus $C$ intersects other components of $D$, then
$H^{0}(\det \mn_{C}\otimes \ml^{-1})$ consists of sections of $\ml$ satisfying
certain additional conditions. For instance if $D=D_{1}+D_{2}+D_{3}$, $D_{1}$
and $D_{2}$ intersect transversely along $C$ and $D_{3}$ intersects
transversely $C$, then we consider divisors
$D_{\epsilon}=f_{1}f_{2}f_{3}+(\epsilon f)^{2}$, where $f\in H^{0}(\ml)$ is
any section vanishing along $f_{1}=f_{2}=f_{3}$. If there are many such
sections then as before the dimension of the singular set goes down. As the
singular set we get  the intersection of $C$ with $f=0$.
Singularities of $D_{\epsilon}$ at points where $f_{1}=f_{2}=f=0,\;f_{3}\not=0$ have the same
type as before ($A_{1}$) but 
at points of $D_{1}\cap D_{1}\cap D_{3}$ we get singularities  of type
$D_{4}$ in general.

Now consider  a triple subvariety  $C$ of the divisor $D$. The transverse
deformations of the double covers correspond to divisors which are also
singular along $C$, so this time the singular set does not decrease but the
type of singularity can change. If $D=D_{1}+D_{2}+D_{3}$, where $D_{i}$ are
smooth divisors such that  $D_{i}$ and $D_{j}$ intersects transversely along $C$ then
$H^{0}(\mn _{C}\otimes\ml^{-1})$ consists of the sections of $\mn_{C}$ that
vanish along $C$. 
If $C$ is a component of the triple point locus  which intersects some components of $D$ that do
not contain $C$, then $H^{0}(\mn _{C}\otimes\ml^{-1})$ contains the sections
of $\mn_{C}$ that vanish along $C$ and satisfy additional vanishing conditions
(of higher order) at the intersection points.

More generally if the multiplicity of $D$ along a subvariety (this time of
arbitrary codimension) is an odd number $2p+1$, then after blowing--up we add
to the branch locus the exceptional divisor. The transverse deformations
corresponds to the divisors that have multiplicity $2p$ along $C$.

In the following series of examples we shall see some of the possible
phenomena occuring for divisors in $\PP$. In higher dimension the situation can be
much more complicated.
\begin{example}
Let $D=D_{1}+D_{2}+D_{3}$ be a sum of three surfaces in $\PP$ of degree
resp. $d_{1},\;d_{2},\;d_{3}$. Let $f_{i}$ be a homogeneous equation of $D_{i}$. 

For $d_{1}=d_{2}=1$, $d_{3}=2$ the intersection $D_{1}\cap D_{2}\cap D_{3}$ contains two
points. For any degree 2 form vanishing at these two points $f$ we consider a
divisor $D_{\epsilon}=f_{1}f_{2}f_{3}+(\epsilon f)^{2}$. For generic choice of
$f$ the divisor $D_{\epsilon}$ has four nodes (the points of intersection of
conics $D_{1}\cap D_{3}$ and $D_{2}\cap D_{3}$ with $f=0$ not lying on the
line $D_{1}\cap D_{2}$) and additional two double points (of type $D_{4}$) at
$D_{1}\cap D_{2}\cap D_{3}$.

 If $d_{1}=d_{2}=1$, $d_{3}=4$ then $D_{1}\cap D_{2}\cap D_{3}$ contains four
points, so every degree 3 form that contains them contains the line $D_{1}\cap
D_{2}$. For a generic choice of such cubic $f$ the divisor
$D_{\epsilon}=f_{1}f_{2}f_{3}+(\epsilon f)^{2}$ has a double line $D_{1}\cap D_{2}$
and 16 nodes (the points of intersection of
quartics $D_{1}\cap D_{3}$ and $D_{2}\cap D_{3}$ with $f=0$ not lying on the
double line).

Assume that $d_{1}=d_{2}=d_{3}=2$ and the forms $f_{i}$ are  dependent. Then $D_{i}$'s
are elements of a pencil of quadrics containing a fixed elliptic curve $C$. For a generic
cubic form $f$ vanishing at $C$ the divisor
$D_{\epsilon}=f_{1}f_{2}f_{3}+(\epsilon f)^{2}$ has double points at $C$
(c-$A_{2}$ singularities). 

Again take $f_{1}$, $f_{2}$ and $f_{3}$ three quadrics containing a smooth
elliptic curve $C$ and let $f_{4}$ be a generic quadric. For a quartic
form $f$
which vanishes along $C$ and has double zeros at the eight points of intersection
of $C$  with $D_{4}$ the divisor $D_{\epsilon}=f_{1}f_{2}f_{3}f_{4}+(\epsilon
f)^{2}$  have eight
fourfold points, it can be written in the form
$G_{4}(f_{1},f_{2},f_{4})$, where $G_{4}$ is a quartic form
(cf. Example~\ref{ex:fourfold}). All  the transverse deformations can be
written as $f_{1}f_{2}f_{3}f_{4}+\epsilon^{2}g$, where $g$ is an octic form
with multiplicity 2 along $C$ and 
eight fourfold points at $D_{1}\cap D_{2}\cap D_{3}\cap D_{4}$. Observe
that the space of such octic forms has dimension 12. It contains the space of
transverse deformations of dimension 7 and a fivedimensional subspace of the
space of equisingular deformations (those octics that can be written as a
degree four polynomial in $f_{1}$ and $f_{2}$).

Take $D_{1}$, $D_{2}$, $D_{3}$  quadrics intersecting at 8 points and
$D_{4}$ a generic quadric vanishing at these points ($f_{4}$ is a linear combination of
$f_{1}$, $f_{2}$, $f_{3}$). For a quartic form $f$
which has double zeros at the eight points of intersection of $f_{i}$'s 
the divisor $D_{\epsilon}=f_{1}f_{2}f_{3}f_{4}+(\epsilon f)^{2}$ has eight
ordinary fourfold points, it can be written in the form
$G_{4}(f_{1},f_{2},f_{3})$, where $G_{4}$ is a quartic form
(cf. Example~\ref{ex:fourfold}). 
In this example transverse deformations can be
written as $f_{1}f_{2}f_{3}f_{4}+\epsilon^{2}g$, where $g$ is an octic form
with eight fourfold points at $D_{1}\cap D_{2}\cap D_{3}\cap D_{4}$.
\end{example}

\begin{remark}\label{rem52}
  The space $H^{1}(\Theta_{\tilde Y}\otimes \tilde\ml^{-1})$ contains deformations
  of $\tilde X$ which are not a double cover of a deformation of $Y$. On
  the other hand if $Y=\PP[n]$ and $D$ is a degree $d$ hypersurface then
  $H^{1}(\Theta_{\PP[n]}(-\frac d2))=0$ (providing $(n,d)\not=(2,6)$).
  From the above description it follows that $H^{1}(\Theta_{\tilde Y}\otimes
  \tilde\ml^{-1})$ also corresponds to deformations of $D$ in $\PP[n]$ but not
  to equisingular ones. So the deformations of $\tilde X$ are smooth models of
  double cover of $\PP[n]$ but not a double cover of a blow--up of
  $\PP[n]$ (cf.\ example~\ref{ex:twosurf}). 
\end{remark}

\section{Deformations of double solids Calabi--Yau threefolds}
\label{sec:cy}

In a special case when $\dim Y=3$, $K_{Y}\cong\ml^{-1}$ and $k=0$ we get 
$h^{1}(\mt_{\tilde Y}\otimes\tilde \ml^{-1})= h^{1}(\tyl)$ if $C$ is a point and 
$h^{1}(\mt_{\tilde Y}\otimes\tilde \ml^{-1})= h^{1}(\tyl)+g(C)$, where $C$ is a
curve ($g(C)$ denotes the genus of $C$).
Now, if we consider an octic surface $D\subset\PP$ and find a resolution of
the double cover induced by a sequence $\sigma:\tilde Y\lra\PP$ of blow--ups of
fourfold and fivefold points and double and triple smooth curves
then $h^{1}(\Theta_{\tilde Y}\otimes\tilde {\mathcal L})$ 
is the sum of genera of all blown-up curves. In this
situation the double cover $\pi:\tilde X\lra\tilde Y$ is a Calabi--Yau manifold. Every
blow--up of a curve gives rise to a ruled surface in $\tilde X$. For a ruled
surface $E\subset \tilde X$ over a genus $g>1$ curve, $E$ deforms with $X$ on a submanifold
of codimension $g$ of the Kuranishi space of $\tilde X$,
over a general point of the Kuranishi space $E$ is replaced by a sum of $2g-2$
rational curves (see \cite{Wilson,sze}). 

\begin{remark}
  Theorem~\ref{thm2} and the above description allow us to compute the number
  of infinitesimal deformations, and consequently also the Hodge numbers of
  the Calabi--Yau manifold $\tilde X$. To compute the
  number of equisingular deformations of the branch locus we can use a computer
  algebra system.
We give two explicit examples, in one the dimension of the space of
equisingular deformations can be computed directly, in  the other, using
Theorem~\ref{thm2} and a Singular program.
\end{remark}

\begin{example}\label{ex:fourfold}
  Let $D=\{(x^{2}-z^{2})^{4}+(y^{2}-w^{2})^{4}+(z^{2}-w^{2})^{4}=0\}$. Then
  $D$ is an irreducible octic with with eighth ordinary fourfold points in
  the vertices of a cube. One easily verify that the space of octics
  with fourfold points in the  vertices of a cube has dimension 14 (and
  there are only finitely many automorphisms of $\PP$ that fixes the vertices of cube,
  namely the symmetries of an affine cube). 

We can deform this octic to another octic with eight 4--fold points if they
are intersection of three quadrics, so we can take generic seven points and
then the eighth is determined. This means that the kernel of the map 
$H^{1}\mt_{\tilde Y}\lra  H^{1}\mn_{D^{\ast}|\tilde Y}$ has dimension $6$ and
\(H^{1}\mt_{\tilde X}\cong H^{1}(\mt_{\tilde Y}(\log D^{\ast}))\cong\CC^{20} \).
Since $\tilde X$ is a Calabi--Yau manifold 
 $H^{1}\mt_{\tilde X}\cong H^{1}\omtx[2]$.
Moreover 
$e(X)=-8$ and so we get $H^{1}\omtx\cong\CC^{16}$. 
Since the group of symmetric (w.r.t natural involution) divisors has rank $9$,
the rank of the skew--symmetric part of the Picard group is $7$.
\end{example}

\begin{example}\label{ex14}
    Let $D$ be the image of generic abelian surface of type $(1,4)$ by the
    mapping defined by the polarization. Surfaces of this type studied in the
    paper \cite{BLS}. The octic  $D$ has four fourfold points and a double
    curve which is a union of four rational curves. The singularities of $D$ can
    be resolved by blowing first the four fourfold points and then the double
    curves (which after the first blow--up are disjoint and smooth). So
    $h^{1}(\mt_{\tilde Y}\otimes \tilde{\ml}^{-1})=0$. 
    
    The equation of $D$ depends on three parameters, we shall
    consider explicit example with
    $\lambda_{0}=\lambda_{1}=\lambda_{2}=\lambda_{3}=1$, the equation takes
    the form 
    \begin{eqnarray*}
      &f=&x^4y^4 +x^4z^4 +z^4t^4+y^4z^4+x^4t^4+y^4t^4
    -2x^2y^2z^4  -2x^4z^2t^2 \\&&+x^2y^2z^2t^2
    +2y^4z^2t^2 +4y^2z^4t^2  +2x^2y^2t^4  -4x^2z^2t^4       
    \end{eqnarray*}
    The fourfold points of $f$ have coordinates
    $(1:0:0:0),(0:1:0:0)$, $(0:0:1:0),(0:0:0:1)$, the double curves are given
    by 
    \begin{eqnarray*}
      &&x=y^2 z^2+y^2 t^2+z^2 t^2=0\\
      &&y=x^2 z^2+x^2 t^2+z^2 t^2=0\\
      &&z=y^2 x^2+y^2 t^2+x^2 t^2=0\\
      &&t=y^2 z^2+y^2 x^2+z^2 x^2=0
    \end{eqnarray*}
Since the double curves have nodes as the only singularities the
symbolic powers coincide of their ideals with usual powers. 
If we consider local coordinates $(x,y,z)$ around the point
$(0:0:0:1)$, then the blow--up at this point is given locally by the
maps $(x,y,z)\mapsto(x,xy,xz)$,
$(x,y,z)\mapsto(xy,x,yz)$ and $(x,y,z)\mapsto(xz,yz,z)$. Taking the
Jacobi matrices of these maps and representing them in the coordinates
on $\PP$ we get the following rational vectorfields $\frac
1x\frac{\partial}{\partial x}$, $\frac
1y\frac{\partial}{\partial y}$, $\frac
1z\frac{\partial}{\partial z}$.  

To compute the dimension of the space of equisingular
deformations we use the following program in Singular

\begin{verbatim}
ring r=0,(x,y,z,t),dp;
poly
octic=x^4*y^4+x^4*z^4-2x^2*y^2*z^4+y^4*z^4-2x^4*z^2*t^2+\
            x^2*y^2*z^2*t^2+2y^4*z^2*t^2+4*y^2*z^4*t^2+x^4*t^4+\
            2x^2*y^2*t^4+y^4*t^4-4*x^2*z^2*t^4+z^4*t^4;
ideal jff=jacob(octic);
ideal jf=ideal(jff[1]/x,jff[2]/y,jff[3]/z,jff[4]/t);
ideal i1=std((x,y,z)^4+jff);
ideal i2=std((x,y,t)^4+jff);
ideal i3=std((x,z,t)^4+jff);
ideal i4=std((y,z,t)^4+jff);
ideal i5=std((ideal(x,y^2*z^2+y^2*t^2+z^2*t^2))^2+\
             ideal(jff[1],jf[2],jf[3],jf[4])); 
ideal i6=std((ideal(y,x^2*z^2+x^2*t^2+z^2*t^2))^2+\
             ideal(jf[1],jff[2],jf[3],jf[4])); 
ideal i7=std((ideal(z,y^2*x^2+y^2*t^2+x^2*t^2))^2+\
             ideal(jf[1],jf[2],jff[3],jf[4])); 
ideal i8=std((ideal(t,y^2*z^2+y^2*x^2+z^2*x^2))^2+\
             ideal(jf[1],jf[2],jf[3],jff[4])); 
ideal ieq=std(intersect (i1,i2,i3,i4,i5,i6,i7,i8)); 
int s=0;
for (int i=1;i<=9;i++)
{s=s+hilb(std(jff),2)[i]-hilb(ieq,2)[i];};
s;
\end{verbatim}
from which we get $h^{1}(\mt_{\tilde Y}(\log D^{\ast}))=3$ and consequently
$h^{1}(\mt_{\tilde X})=h^{1,2}(\tilde X)=3$. 
As the Euler characteristic of $\tilde X$ is easily computed to be 24 we
have $h^{1,1}(X)=\rho(X)=15$. It is easy to see that the group of
symmetric divisors on $X$ has rank 9 (pullback of a plane in $\PP$ and 8
exceptional divisors of blow--ups), so the rank of the group of
skew-symmetric divisors~is~6.

The dimension of transverse deformations is 3, we get the same result
if we do not consider the rational vectorfields in the formula for
transverse deformations. The explanation is that any component of the
double points locus is a rational quartic with four nodes. After
blow--up of nodes we get rational curves, which after deformation and
projecting to $\PP$ are quartic with three nodes. Since the quartic
with three nodes is uniquely determined by the nodes and the tangent
lines at nodes, the deformation can be realized as a deformation in $\PP$.
\end{example}

\subsection*{Acknowledgement}
The work was done during the first named author stay at the Johannes
  Gutenberg--Universit\"at in Mainz supported by the European
  Commission grant no HPRN-CT-2000-00099 and partially by DFG
  Schwerpunktprogramm 1094 (Globale Methoden in der komplexen Geometrie).

\end{document}